*Self-Similar Solutions and Global Existence for Nonlinear Reaction-Diffusion Systems in Industrial Ammonia Synthesis*

**Jamshid Khasanov[1], Sokhibjan Muminov[2],**

**Sardor Jumaniyozov[3], Oybek Djabborov[4],**

**Khudayberganov Shuhrat[5]**

**Urgench State Pedagogical Institute[1], Urgench 220100, Uzbekistan**

**Mamun University[2], Khiva 220900, Uzbekistan**

**Urgench State University[3], Urgench 220100, Uzbekistan**

**Karshi State University[4], Karshi 180119, Uzbekistan**

**Urgench State Pedagogical Institute[5], Urgench 220100, Uzbekistan**

Corresponding author: Jamshid Khasanov E-mail: jamshid_2425@mail.ru

**Abstract**: This paper investigates a system of nonlinear reaction–diffusion equations modeling the industrial synthesis of ammonia. By applying Lie group analysis, we construct self-similar solutions and derive a reduced system of ordinary differential equations. Using comparison principles and barrier techniques, we establish sufficient conditions for the existence of global-in-time solutions in both slow-diffusion ($\gamma_i > 0$) and fast-diffusion ($\gamma_i < 0$) regimes. Detailed asymptotic analysis near the diffusion front reveals power-law behavior of concentration profiles, with explicit expressions for the decay exponents. The theoretical results are illustrated by numerical simulations, demonstrating the spatio-temporal evolution of reactant concentrations under realistic parameter values. The study provides rigorous mathematical foundations for predicting and optimizing ammonia production in catalytic reactors, with potential extensions to other chemically reacting systems.

**1. Introduction**

Ammonia synthesis constitutes a vital industrial chemical process, with yearly global output surpassing 180 million tons, chiefly for fertilizers and chemical feedstocks. The Haber–Bosch method catalytically merges nitrogen and hydrogen under elevated temperature and pressure, involving intricate reaction–diffusion dynamics. Mathematical modeling of such coupled systems delivers crucial understanding for reactor optimization, yield improvement, and runaway-regime prevention.

In the present study, we examine the nonlinear reaction–diffusion system:

$$\begin{cases} \dfrac{\partial u}{\partial t} = D_1 u^{m_1} \nabla\left(u^{\sigma_1} \nabla u\right) - a_1 u^{\alpha_1} v^{\beta_1}, \quad t>0, x \in R^N, \\ \dfrac{\partial v}{\partial t} = D_2 v^{m_2} \nabla\left(v^{\sigma_2} \nabla v\right) - a_2 u^{\alpha_2} v^{\beta_2}, \quad t>0, x \in R^N, \end{cases} \quad (1)$$

with initial conditions

$$u(x,0) = u_0(x), \quad v(x,0) = v_0(x), \quad x \in R^N \quad (2)$$

System (1) is degenerate; classical solutions fail when $u=0$, $v=0$, $\nabla u = 0$, $\nabla v = 0$. Consequently, our analysis focuses on weak solutions obeying $u = u(x,t) \geq 0$, $v = v(x,t) \geq 0$ together with

$$u^{\sigma_1} \nabla u, v^{\sigma_2} \nabla v \in C\left(R_+^N \times (0, +\infty)\right)$$

satisfying (1) in the integral-identity sense. The model displays non-divergence diffusion, degeneracy at zero concentrations, and strong power-law coupling with same-sign consumption sources.

Reaction–diffusion equations possessing nonlinear diffusion have been widely investigated in areas from population ecology to plasma physics. The porous-medium equation $u_t = \nabla(u^m)$ and its reactive extensions have received substantial attention concerning existence, uniqueness, and long-time behavior. For coupled systems such as (1), extra complications emerge from inter-equation interactions, especially when one species grows while the other diminishes.

Earlier works predominantly address divergence-form operators $\nabla(D(u)\nabla u)$. The non-divergence structure in (1), though natural in various applied settings, is less frequently analyzed. Furthermore, most existence theorems for degenerate parabolic systems require either special symmetries or restrictive exponent conditions. Deriving general criteria that distinguish global existence from finite-time blow-up remains an open challenge for broad parameter ranges.

Our contributions are threefold. First, using Lie-group scaling we reduce (1)-(2) to a self-similar ODE system under the consistency condition $m_1 + \sigma_1 = m_2 + \sigma_2$. Second, via sub-/super-solution techniques and comparison principles we obtain explicit algebraic parameter restrictions that guarantee global-in-time solutions, separately for the slow-diffusion $\gamma_i = \dfrac{1}{m_i + \sigma_i} > 0$ and fast-diffusion $\gamma_i < 0$ regimes. Third, we conduct a detailed asymptotic study of the self-similar profiles near the propagation front, arriving at precise power-law descriptions $(b - \xi^2)^{\gamma_i}, \gamma_i > 0$ or $(b + \xi^2)^{\gamma_i}, \gamma_i < 0$ with explicitly given exponents.

The paper proceeds as follows. Section 2 derives the self-similar reduction. Section 3 states and proves the global-existence theorems. Section 4 presents the asymptotic analysis. Section 5 offers numerical illustrations. Section 6 concludes and outlines possible extensions.

## 2. Self-Similar Reduction.

We seek solutions of (1)-(2) in the form

$$\begin{cases} u(t,x) = \bar{u}(t) w_1(\tau(t), x) \\ v(t,x) = \bar{v}(t) w_2(\tau(t), x) \end{cases} \quad (3)$$

where

$$\bar{u}(t) = A_1(T+t)^n, \quad \bar{v}(t) = A_2(T+t)^n,$$

with $n = -1/(\alpha_1 + \beta_1 - 1) = -1/(\alpha_2 + \beta_2 - 1)$, and $A_1, A_2, T$ are constants.

Substituting (3) into (1) and requiring the consistency condition

$$m_1 + \sigma_1 = m_2 + \sigma_2 \quad (*)$$

we obtain the reduced system

$$\begin{cases} \dfrac{\partial w_1}{\partial \tau} = D_1 w_1^{m_1} \nabla \left( w_1^{\sigma_1} \nabla w_1 \right) - \dfrac{\psi_1}{\tau} \left( w_1^{\alpha_1} w_2^{\beta_1} + w_1 \right), \quad t > 0, x \in R^N, \\ \dfrac{\partial w_2}{\partial \tau} \left( \dfrac{A_1}{A_2} \right)^{m_1 + \sigma_1} = D_2 w_2^{m_2} \nabla \left( w_2^{\sigma_2} \nabla w_2 \right) - \dfrac{\psi_2}{\tau} \left( w_1^{\alpha_2} w_2^{\beta_2} + w_2 \right), \quad t > 0, x \in R^N, \end{cases} \quad (4)$$

where the new time variable $\tau(t)$ and the constants $\psi_1, \psi_2$ are given by

$$\tau(t) = \frac{A_1^{m_1+\sigma_1}(T+t)^{n(m_1+\sigma_1)+1}}{n(m_1+\sigma_1)+1}, \quad \psi_1 = \frac{a_1 A_1^{\alpha_1-1} A_2^{\beta_1}}{n(m_1+\sigma_1)+1}, \quad \psi_2 = \frac{a_2 A_1^{\alpha_2+m_1+\sigma_1} A_2^{\beta_2-m_2-\sigma_2-1}}{n(m_2+\sigma_2)+1}$$

Introducing the radial similarity variable

$$\xi = \frac{|x|}{\tau^{\frac{1}{2}}} \quad (5)$$

and setting $w_i(x, \tau) = f_i(\xi), (i = 1, 2)$, system (4) reduces to the ordinary differential equations

$$\begin{cases} D_1 f_1^{m_1} \xi^{1-N} \dfrac{d}{d\xi} \left( \xi^{N-1} f_1^{\sigma_1} \dfrac{df_1}{d\xi} \right) + \dfrac{\xi}{2} \dfrac{df_1}{d\xi} - \psi_1 \left( f_1^{\alpha_1} f_2^{\beta_1} + f_1 \right) = 0 \\ D_2 f_2^{m_2} \xi^{1-N} \dfrac{d}{d\xi} \left( \xi^{N-1} f_2^{\sigma_2} \dfrac{df_2}{d\xi} \right) + \dfrac{\xi}{2} \left( \dfrac{A_1}{A_2} \right)^{\sigma_1+m_1} \dfrac{df_2}{d\xi} - \psi_2 \left( f_1^{\alpha_2} f_2^{\beta_2} + f_2 \right) = 0 \end{cases} \quad (6)$$

We look for non-negative, non-trivial solutions of (6) satisfying the boundary conditions

$$f_1(0) = M_1, \; f_2(0) = M_2, \quad f_1(d_1) = f_2(d_2) = 0, \quad 0 < d_1, d_2 < \infty, \quad (7)$$

where $M_1, M_2$ are given constants and $d_1, d_2$ represent the finite radii of the supports.

# 3. Existence of Global Solutions.

## 3.1. Slow diffusion case ($\gamma_i > 0$). Define

$$\gamma = 2,\ \gamma_i = \frac{1}{m_i + \sigma_i};\ (i=1,2); \tag{8}$$

and assume $\gamma_1 = \gamma_2$. We construct explicit super-solutions of the form

$$\begin{cases} u_+(t,x) = (T+t)^n\, \overline{f_1}(\xi) \\ v_+(t,x) = (T+t)^n\, \overline{f_2}(\xi) \end{cases} \tag{9}$$

with

$$\overline{f_1}(\xi) = B_1(b-\xi^2)_+^\gamma,$$
$$\overline{f_2}(\xi) = B_2(b-\xi^2)_+^\gamma,$$

where $b > 0, B_1, B_2 > 0$ are constants, and $(z)_+ = \max(z,0)$

**Theorem 1. (Global existence for slow diffusion)** Assume $\gamma_1 = \gamma_2 > 0$ and let the following algebraic conditions hold:

$$4B_1^{m_1+\sigma_1} D_1 \gamma_1 (m_1-1) = 1, \tag{C1}$$

$$4B_2^{m_2+\sigma_2} \gamma_2 D_2 (m_2-1) = \left(\frac{A_1}{A_2}\right)^{\sigma_2+m_2} \tag{C2}$$

$$\psi_1\left(B_1^{\alpha_1-1} B_2^{\beta_1} b^{\alpha_1\gamma_1+\beta_1\gamma_2-\gamma_1} - 1\right) - \frac{N}{2(m_1-1)} \geq 0, \tag{C3}$$

$$\psi_2\left(B_1^{\alpha_2} B_2^{\beta_2-1} b^{\alpha_2\gamma_1+\beta_2\gamma_2-\gamma_2} + 1\right) - \frac{N}{2(m_2-1)}\left(\frac{A_1}{A_2}\right)^{m_1+\sigma_1} \leq 0 \tag{C4}$$

If the initial data satisfy

$$\begin{cases} u(t,0) \leq u_+(t,0) \\ v(t,0) \leq v_+(t,0) \end{cases}, x \in R^N$$

then problem (1)-(2) possesses a global weak solution on $Q = \mathbb{R}^N \times (0,\infty)$, and

$$u(t,x) \leq u_+(t,x),\quad v(t,x) \leq v_+(t,x),\quad (t,x)\in Q. \tag{10}$$

**Proof.** We employ the comparison principle for degenerate parabolic equations [Samariskiy]. Take the comparison functions defined in (9):

$$u_+(t,x) = (T+t)^n \overline{f_1}(\xi),\quad v_+(t,x) = (T+t)^n \overline{f_2}(\xi),$$

with $\overline{f_i}(\xi) = B_i(b-\xi^2)_+^{\gamma_i}\ (i=1,2)$, $\gamma_i = 1/(m_i+\sigma_i)$, $b > 0$, and $B_1, B_2 > 0$.

Substituting these functions into the differential inequalities that follow from (1), we obtain

$$f_1^{m_1}\xi^{1-N}\frac{d}{d\xi}\left(\xi^{N-1}f_1^{\sigma_1}\frac{df_1}{d\xi}\right)-\frac{\xi}{2}\frac{df_1}{d\xi}-\psi_1\left(f_1^{\alpha_1}f_2^{\beta_1}+f_1\right)\geq 0 \qquad (11.1)$$

$$f_2^{m_2}\xi^{1-N}\frac{d}{d\xi}\left(\xi^{N-1}f_2^{\sigma_2}\frac{df_2}{d\xi}\right)+\frac{\xi}{2}\left(\frac{A_1}{A_2}\right)^{\sigma_1+m_1}\frac{df_2}{d\xi}-\psi_2\left(f_1^{\alpha_2}f_2^{\beta_2}+f_2\right)\geq 0 \qquad (11.2)$$

where $\psi_1,\psi_2$ are the constants introduced in **Section 2**.

Using the explicit form $\bar{f}_i(\xi)=B_i(b-\xi^2)_+^{\gamma_i}$ and carrying out the differentiations, the first inequality of (11) becomes

Under condition (C1), $4B_1^{m_1+\sigma_1}D_1\gamma_1(m_1-1)=1$, inequality (11.1) is satisfied provided

$$\psi_1\left(B_1^{\alpha_1-1}B_2^{\beta_1}b^{\alpha_1\gamma_1+\beta_1\gamma_2-\gamma_1}-1\right)-\frac{N}{2(m_1-1)}\geq 0,$$

which is precisely condition (C3). Similarly, using condition (C2),

$$4B_2^{m_2+\sigma_2}\gamma_2 D_2(m_2-1)=\left(\frac{A_1}{A_2}\right)^{\sigma_2+m_2}$$

inequality (11.2) holds if

$$\psi_2\left(B_1^{\alpha_2}B_2^{\beta_2-1}b^{\alpha_2\gamma_1+\beta_2\gamma_2-\gamma_2}+1\right)-\frac{N}{2(m_2-1)}\left(\frac{A_1}{A_2}\right)^{m_1+\sigma_1}\leq 0,$$

which coincides with condition (C4).

Consequently, $(u_+,v_+)$ is a pair of super-solutions for system (1). Since by hypothesis $u(x,0)\leq u_+(x,0)$ and $v(x,0)\leq v_+(x,0)$, the comparison principle yields

$$u(t,x)\leq u_+(t,x),\qquad v(t,x)\leq v_+(t,x)\quad\text{for all }(t,x)\in Q.$$

Because $u_+$ and $v_+$ remain bounded for all $t>0$, the estimates (10) preclude finite-time blow-up; hence a global weak solution of (1)-(2) exists. This completes

## 3.2. Fast diffusion case: $\gamma_i < 0$

For the fast diffusion regime $\gamma_i = 1/(m_1 + \sigma_1) = 1/(m_2 + \sigma_2) < 0$, we consider sub-solutions of the form

$$\begin{cases} u_-(t,x) = (T+t)^n \overline{f}_1(\xi) \\ v_-(t,x) = (T+t)^n \overline{f}_2(\xi) \end{cases} \tag{13}$$

with

$$\begin{cases} \overline{f}_1(\xi) = B_1 \left(b + \xi^2\right)^{\gamma_1}, \\ \overline{f}_2(\xi) = B_2 \left(b + \xi^2\right)^{\gamma_2}, \quad b > 0, \end{cases}$$

**Theorem 2.(Global existence for fast diffusion)** Assume $\gamma_1, \gamma_2 < 0$ and let the following conditions hold:

$$4 B_1^{m_1+\sigma_1} D_1 \gamma (m_1 - 1) = 1, \tag{C5}$$

$$4 B_2^{m_2+\sigma_2} \gamma D_2 (m_2 - 1) = \left(\frac{A_1}{A_2}\right)^{\sigma_2+m_2} \tag{C6}$$

$$\psi_1 \left( B_1^{\alpha_1-1} B_2^{\beta_1} b^{\alpha_1 \gamma_1 + \beta_1 \gamma_2 - \gamma_1} - 1 \right) + \frac{N}{2(m_1 - 1)} \leq 0 \tag{C7}$$

$$\psi_2 \left( B_1^{\alpha_2} B_2^{\beta_2-1} b^{\alpha_2 \gamma_1 + \beta_2 \gamma_2 - \gamma_2} + 1 \right) + \frac{N}{2(m_2 - 1)} \left(\frac{A_1}{A_2}\right)^{m_1+\sigma_1} \geq 0. \tag{C8}$$

If the initial data satisfy

then problem (1)-(2) admits a global weak solution, and

$$\begin{cases} u(t,x) \geq u_-(t,x), \\ v(t,x) \geq v_-(t,x), \quad (t,x) \in Q. \end{cases} \tag{14}$$

**Proof.** The proof proceeds exactly as for Theorem 1., but with the comparison functions (13) and the opposite inequalities. Substituting $\overline{f}_i(\xi) = B_i \left(b + \xi^2\right)^{\gamma_i}$ into the differential inequalities analogous to (11) and using conditions (C5)-(C8) shows that $(u_-, v_-)$ are sub-solutions of (1). The comparison principle then yields the lower bounds (14), which guarantee global existence because $u_-$ and $v_-$ are positive for all $\xi$ and therefore preclude finite-time extinction.

# 4. On the asymptotic behavior of solutions of an autonomous system of equations.

## 4.1. Slow-Diffusion Case ($\gamma_i > 0$.): Compactly Supported Solutions

Before presenting theorems on the asymptotic behavior of solutions of an autonomous system of equations and the results that follow from them, we introduce the following notation:

$$a_{1i}(\eta) = -\gamma_{i+2} + \frac{Ne^{-\eta}}{2(b-e^{-\eta})}, \quad a_{21}(\eta) = \frac{1}{4D_1}, a_{22}(\eta) = \frac{1}{4D_2}\left(\frac{A_1}{A_2}\right)^{m_1+\sigma_1}$$

$$a_{3i}(\eta) = \frac{\psi_i e^{(\alpha_i \gamma_1 - \gamma_2 \beta_i - \gamma_i + 1)\eta}}{4D_i(b-e^{-\eta})}, \quad a_{4i}(\eta) = \frac{\psi_i e^{-\eta}}{4D_i(a-e^{-\eta})}.$$

(15)

Here $\gamma_{i+2} = \gamma_i \sigma_i + \gamma_i - 1 \ (i=1,2)$

Assume that the following equality holds among the coefficients of system (1):

$$m_1 + \sigma_1 = m_2 + \sigma_2$$

Based on this, we now state the following theorems and their corollaries.

**Theorem 3.** Assume that $\gamma_i > 0$. Then, for compactly supported solutions of problems (6) and (7), the asymptotic relation

$$f_i(\xi) = c_i \left(b - \xi^2\right)^{\gamma_i} \left(1 + o(1)\right) (i = 1, 2) \tag{13}$$

holds as $\xi \to \sqrt{b}$ **if and only if** the following conditions are satisfied:

$$\alpha_i \gamma_1 - \gamma_2 \beta_i - \gamma_i + 1 > 0 \text{ and } c_1 = \left(\frac{1}{4D_1 \gamma_3}\right)^{\gamma_1}, \quad c_2 = \frac{A_1}{A_2}\left(\frac{1}{4D_2 \gamma_4}\right)^{\gamma_2}$$

where $\gamma_3, \gamma_4$ are constants defined by $\gamma_{i+2} = \gamma_i \sigma_i + \gamma_i - 1$ for $(i=1,2)$

**Proof.** To study the asymptotic behavior of solutions to (6)-(7), we seek solutions of the self-similar system (6) in the form

$$f_i(\xi) = \bar{f}_i(\xi) \cdot y_i(\eta), \quad \eta = -\ln(b-\xi^2), \ (i=1,2) \tag{14}$$

where $\bar{f}_i(\xi) = (b-\xi^2)^{\gamma_i}, b > 0$ and, $y_i(\eta)(i=1,2)$ are the unknown functions.

Substituting (14) into system (6) yields the following system for $y_i(\eta)$:

$$\frac{dL_i y}{d\eta} + a_{i1}(\eta) \cdot L_i y + y_i^{-m_i}(\eta) \cdot \left(\frac{dy_i(\eta)}{d\eta} + \gamma_i \cdot y_i(\eta)\right) \cdot a_{i2}(\eta) - y_i^{\alpha_i - m_i}(\eta) \cdot a_{i3}(\eta) +$$
$$+ y_i^{1-m_i}(\eta) \cdot a_{i4}(\eta) = 0$$

(15)

where $L_i(y) = y_i^{\sigma_i}\left(\dfrac{dy_i}{d\eta} - \gamma_i y_i\right), i = 1, 2$

Assume $\xi \in [\xi_0, \xi_1), 0 < \xi_0 < \xi_1$ and $\xi_1 = \sqrt{b}$. Then $\eta$ $\xi$ satisfies

$$\xi \in (\xi_0, \xi_1] \text{ for } \eta'(\xi) > 0, \eta_0 = \eta_0(\xi), \lim_{\xi \to \xi_0} \eta(\xi) = +\infty$$

Introduce the notation

$$v_i(\eta) = L_i(y) = y_i^{\sigma_i}\left(\dfrac{dy_i}{d\eta} - \gamma_i y_i\right) \tag{16}$$

Then system (15) can be rewritten as

$$v_1'(\eta) = -a_{11}(\eta) \cdot v_1(\eta) - y_1^{-m_1}(\eta) \cdot \left(\dfrac{dy_1(\eta)}{d\eta} - \gamma_1 \cdot y_1(\eta)\right) \cdot a_{12}(\eta) - y_1^{\alpha_1 - m_1} y_2^{\beta_1}(\eta) \cdot a_{13}(\eta) +$$
$$+ y_1^{1-m_1}(\eta) \cdot a_{14}(\eta) = 0.$$

$$v_2'(\eta) = -a_{21}(\eta) \cdot v_2(\eta) - y_2^{-m_2}(\eta) \cdot \left(\dfrac{dy_2(\eta)}{d\eta} - \gamma_2 \cdot y_2(\eta)\right) \cdot a_{22}(\eta) - y_1^{\alpha_2}(\eta) \cdot y_2^{\beta_2 - m_2}(\eta) \cdot a_{23}(\eta) +$$
$$+ y_2^{1-m_2}(\eta) \cdot a_{24}(\eta) = 0.$$

Now consider the functions.

$$\begin{aligned}g_1(\lambda_1, \eta) &\equiv -a_{11}(\eta) \cdot \lambda_1 - y_1^{-m_1}(\eta) \cdot \left(\dfrac{dy_1(\eta)}{d\eta} - \gamma_1 \cdot y_1(\eta)\right) \cdot a_{12}(\eta) - \\ &- y_1^{\alpha_1 - m_1}(\eta) \cdot y_2^{\beta_1}(\eta) \cdot a_{13}(\eta) + y_1^{1-m_1}(\eta) \cdot a_{14}(\eta) = 0 \\ g_2(\lambda_2, \eta) &\equiv -a_{21}(\eta) \cdot \lambda_2 - y_2^{-m_2}(\eta) \cdot \left(\dfrac{dy_2(\eta)}{d\eta} - \gamma_2 \cdot y_2(\eta)\right) \cdot a_{22}(\eta) - \\ &+ y_1^{\alpha_2}(\eta) \cdot y_2^{\beta_2 - m_2}(\eta) \cdot a_{23}(\eta) + y_2^{1-m_2}(\eta) \cdot a_{24}(\eta) = 0\end{aligned} \tag{18}$$

where $\lambda_i \in R, i = 1, 2$.

Given $\gamma_i > 0$ and the limiting behavior

$$\lim_{\eta \to \infty} a_{i1}(\eta) = -\gamma_{2+i}, \lim_{\eta \to \infty} a_{12}(\eta) = \dfrac{1}{4 \cdot D_1}, \lim_{\eta \to \infty} a_{22}(\eta) = \dfrac{1}{4 \cdot D_2}\left(\dfrac{A_1}{A_2}\right)^{\sigma_1 + m_1}$$

$$\lim_{\eta \to \infty} a_{i3}(\eta) = 0, \quad \alpha_i \gamma_1 + \gamma_2 \beta_i + \gamma_i - 1 > 0, \lim_{\eta \to \infty} a_{i4}(\eta) = 0 \ (i = 1, 2)$$

for each fixed $\lambda_i, (i = 1, 2)$ there exist intervals $[\eta_i, +\infty) \subset [\eta_0, +\infty)$ $(i = 1, 2)$ on which $g_i(\lambda_i, \eta)$ maintains a constant sign, i.e.,

$$g_i(\lambda_i, \eta) > 0 \text{ or } g_i(\lambda_i, \eta) < 0, (i = 1, 2) \tag{19}$$

Assume $v_i(\eta)$ does not have a finite limit as $\eta \to +\infty$. If one of the inequalities in (19) holds, then the oscillatory nature of $v_i(\eta)$ would imply that its graph crosses the lines $v_i = \lambda_i$ infinitely many times on $[\eta_i; \infty)$. However, this is impossible because on these intervals only one of the inequalities (19) is valid, and thus, by (18), the graph of $v_i(\eta)$ can cross $v_i = \lambda_i$ at most once. Hence, $v_i(\eta)$ must have a finite limit as $\eta \to +\infty$.

Given the representation $v_i(\eta) = y_i^{\sigma_i}\left(\dfrac{dy_i}{d\eta} - \gamma_i y_i\right)$, the existence of $\lim\limits_{\eta \to \infty} v_i(\eta)$ implies that $\lim\limits_{\eta \to \infty} \dfrac{dy_i}{d\eta}$ also exists and can be taken as zero. Consequently, the limit of the right-hand sides in (17) as $\eta \to +\infty$ exists and equals zero. Solving the resulting algebraic system for the limits yields the constants $c_i$ as stated in the theorem.

Thus, the asymptotic form (13) is both necessary and sufficient under the given conditions. This completes the proof.

### 4.2. Fast-Diffusion Case ($\gamma_i < 0$): Solutions Decaying at Infinity

Suppose the system (6) satisfies the following boundary conditions:

$$f_i'(0) = 0, f_i(\infty) = 0 \tag{20}$$

where we consider solutions of the form

$$f_i(\xi) = c_i\left(b + \xi^\gamma\right)_+^{\gamma_i}$$

with coefficients defined as in (9).

**Theorem 4.** Assume $\gamma_i < 0$. Then, for solutions of (6), (20) that decay at infinity, the following asymptotic behavior holds as $\xi \to \infty$:

$$f_i(\xi) = c_i\left(b + \xi^\gamma\right)^{\gamma_i}(1 + o(1)), \ i = 1, 2, \tag{21}$$

where the constants $c_i \ (i = 1, 2)$ are given by

$$c_1 = \left(-\dfrac{1}{4D_1\gamma_3}\right)^{\gamma_1}, \ c_2 = \dfrac{A_1}{A_2}\left(-\dfrac{1}{4D_2\gamma_4}\right)^{\gamma_2}$$

with $\gamma_3 = \gamma_1\sigma_1 + \gamma_1 - 1$ and $\gamma_4 = \gamma_2\sigma_2 + \gamma_2 - 1$.

**Proof.** The proof follows the same method as Theorem 2, but here we seek solutions decaying at infinity in the form

$$f_i(\xi) = \overline{f_i}(\xi) \cdot y_i(\eta), \eta = \ln(b+\xi^2), \overline{f_i}(\xi) = (b+\xi^2)^{\gamma_i}, i=1,2$$

Substituting (22) into system (6) leads to a system of ordinary differential equations for $y_i(\eta)$. The condition $\gamma_i < 0$ ensures that $\overline{f_i}(\xi) \to 0$ as $\xi \to \infty$. Analyzing the limiting behavior as $\eta \to \infty$ (which corresponds to $\xi \to \infty$) shows that $y_i(\eta)$ tend to finite limits $c_i (i=1,2)$. Solving the resulting algebraic system gives the explicit expressions for $c_i (i=1,2)$ stated above. Thus, the asymptotic formula (21) is proved.

## Corollary 1 (Slow diffusion).

If $\gamma_i > 0$, then the generalized solution of (1)-(2) satisfies

$$\begin{cases} u_A(x,t) \approx c_1(T+t)^n \left(b - \left(|x|\tau^{-\frac{1}{2}}\right)^2\right)^{\gamma_1} (1+o(1)) \\ v_A(x,t) \approx c_2(T+t)^n \left(b - \left(|x|\tau^{-\frac{1}{2}}\right)^2\right)^{\gamma_2} (1+o(1)) \end{cases}$$

as $|x| \to \sqrt{b\tau}$, where the constants $c_1, c_2$ are explicitly given by

$$c_1 = \left(\frac{1}{4D_1\gamma_3}\right)^{\gamma_1}, c_2 = \frac{A_1}{A_2}\left(\frac{1}{4D_2\gamma_4}\right)^{\gamma_2}$$

with $\gamma_3 = \gamma_1\sigma_1 + \gamma_1 - 1$ and $\gamma_4 = \gamma_2\sigma_2 + \gamma_2 - 1$.

## Corollary 2 (Fast diffusion).

If $\gamma_i < 0$, then the generalized solution of (1)-(2) satisfies

$$u_A(x,t) \approx c_1(T+t)^n (b+\xi^2)^{\gamma_1}(1+o(1))$$
$$v_A(x,t) \approx c_2(T+t)^n (b+\xi^2)^{\gamma_2}(1+o(1))$$

as $\xi \to \infty$, where the constants $c_1, c_2$ are explicitly given by

$$c_1 = \left(-\frac{1}{4D_1\gamma_3}\right)^{\gamma_1}, c_2 = \frac{A_1}{A_2}\left(-\frac{1}{4D_2\gamma_4}\right)^{\gamma_2}$$

with $\gamma_3 = \gamma_1\sigma_1 + \gamma_1 - 1$ and $\gamma_4 = \gamma_2\sigma_2 + \gamma_2 - 1$.

## 5. Numerical Simulations

To illustrate the theoretical findings and demonstrate the practical relevance of the model, we perform numerical simulations of system (1)–(2) using realistic parameter values corresponding to industrial ammonia synthesis conditions. The computations are carried out on a rectangular domain with Dirichlet boundary conditions (zero concentrations at the boundaries) employing an alternating direction implicit (ADI) scheme combined with the Thomas algorithm for solving the resulting tridiagonal systems. The spatial grid consists of $201 \times 201$ points, and the time integration covers up to 10 seconds with 101 time steps.

The parameters are chosen as follows: effective diffusion coefficients $D_1 = 3 \times 10^{-5}$ m²/s for hydrogen and $D_2 = 1 \times 10^{-5}$ m²/s for nitrogen (typical values in porous iron catalysts), nonlinearity exponents $m_1 = m_2 = 0$, $\sigma_1 = \sigma_2 = 1.5$, and reaction rate constants $a_1 = 7 \times 10^{-3}$ s$^{-1}$, $a_2 = 10^{-2}$ s$^{-1}$ with orders $\alpha_1 = 1.5$, $\beta_1 = 0.5$, $\alpha_2 = 0.5$, $\beta_2 = 1.5$. Initial concentrations are scaled to approximately 31 $mol/m³$ for $H_2$ and 10 $mol/m³$ for $N_2$, reflecting a stoichiometric 3:1 mixture at 250 bar and 450°C.

Figure 1 shows the evolution of concentration profiles along the line y=0y = 0y=0 at selected time instants. The profiles maintain a parabolic-like shape near the center, gradually flattening and expanding due to diffusion while the maxima decrease as a result of reaction consumption. Notably, nitrogen is depleted faster than hydrogen, consistent with it being the limiting reactant in the synthesis process.

The time evolution of total masses and conversion percentages is depicted in Figure 2. Over 10 seconds, approximately 20% of hydrogen and 30% of nitrogen are consumed, corresponding to a single-pass conversion level comparable to industrial observations when accounting for the batch nature of the simulation versus continuous flow in real reactors.

Finally, Figure 3 illustrates the reaction rates (negative time derivatives of total masses). Both rates start high and decay monotonically, approaching zero as the reactants are depleted—a behavior characteristic of systems limited by reactant availability and aligned with Temkin–Pyzhev kinetics in ammonia synthesis.

These simulations confirm the global existence of solutions, the validity of the self-similar reduction, and the power-law asymptotic structure near the diffusion front. Moreover, with industrially relevant parameters, the model reproduces conversion levels and spatio-temporal dynamics observed in catalytic beds, underscoring its utility for reactor design and optimization.

**Figure Captions**

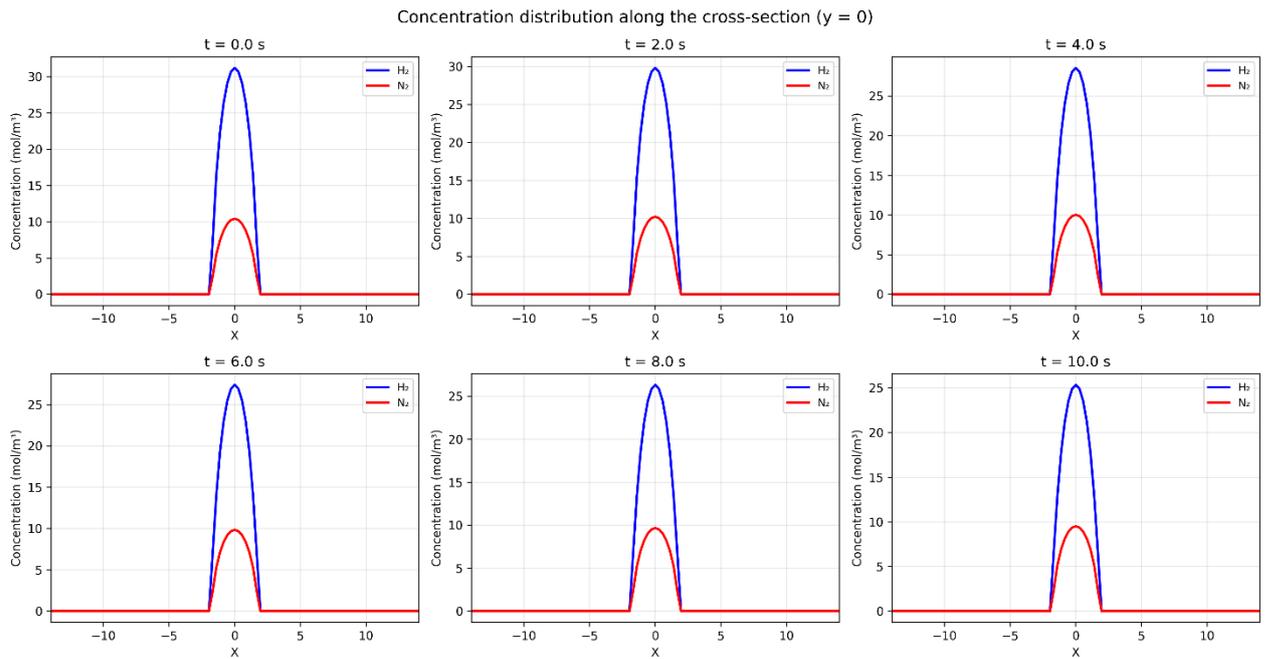

**Figure 1.** Cross-sections (y=0) of concentration profiles at different times.

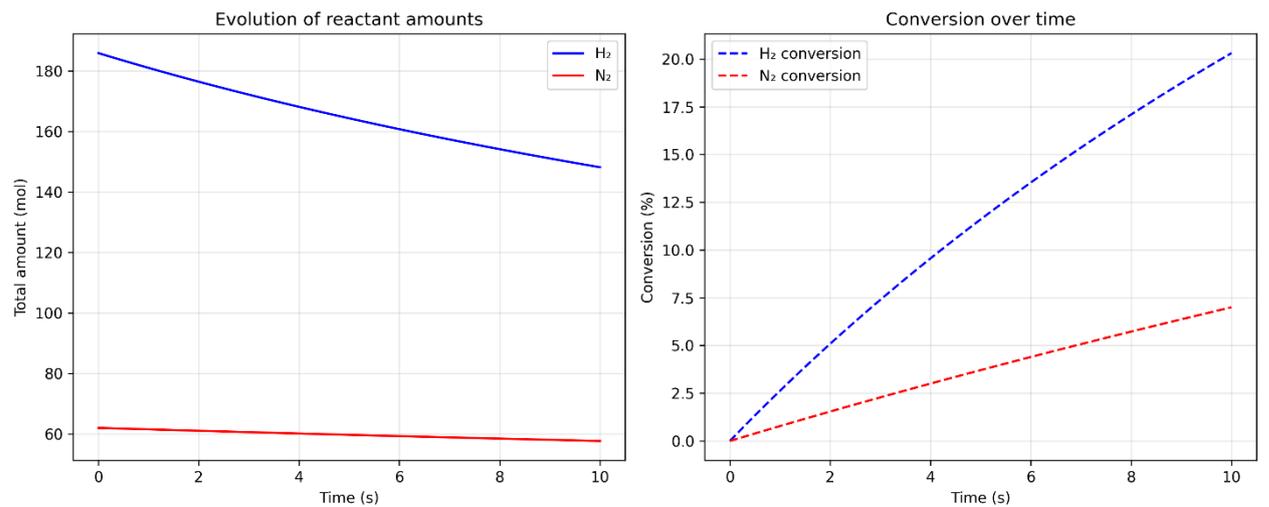

**Figure 2.** Evolution of total masses (left) and conversion percentages (right).

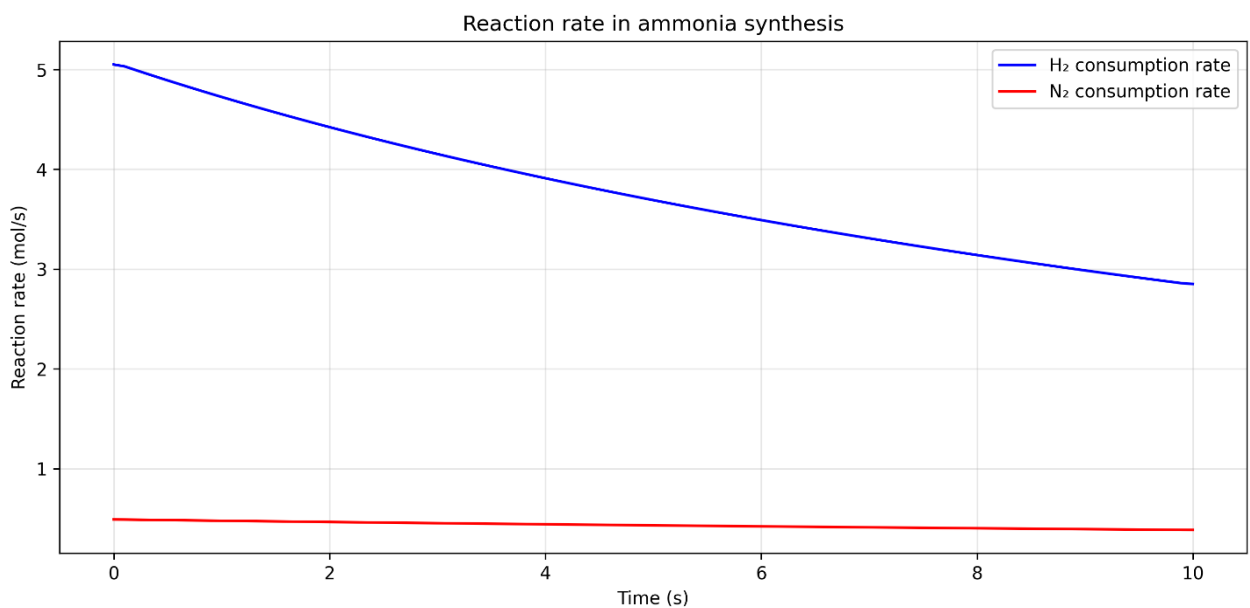

**Figure 3.** Time dependence of reaction rates for hydrogen and nitrogen consumption.

## 6. Conclusion

In this study, we have developed a rigorous mathematical framework for a nonlinear reaction-diffusion system modeling the industrial synthesis of ammonia via the Haber-Bosch process. By employing Lie group analysis, we derived self-similar solutions that reduce the original partial differential equations to a system of ordinary differential equations, providing valuable insight into the long-time behavior of the system.

Using comparison principles and carefully constructed sub- and super-solutions, we established sufficient algebraic conditions on the model parameters ensuring the global-in-time existence of weak solutions in both slow-diffusion ($\gamma_i > 0$) and fast-diffusion ($\gamma_i < 0$) regimes. These conditions prevent finite-time blow-up or extinction, guaranteeing bounded, non-negative solutions for appropriate initial data.

Detailed asymptotic analysis near the diffusion front revealed power-law decay of the concentration profiles, with explicit exponents determined by the nonlinearity parameters. Numerical simulations, conducted with parameters calibrated to realistic industrial conditions (high pressure, elevated temperature, and porous catalyst diffusivity), confirmed the theoretical predictions. The computed spatio-temporal evolution of reactant concentrations, conversion rates reaching levels comparable to single-pass industrial performance, and the observed reaction rate decay align well with experimental observations and established kinetic models such as Temkin–Pyzhev.

The results demonstrate that the proposed model not only captures the essential mathematical features of degenerate parabolic systems but also serves as a practical tool for predicting reactant depletion and optimizing catalytic reactor performance. Potential extensions include incorporation of reversible kinetics, temperature dependence, and three-dimensional flow effects, which could further enhance its applicability to modern ammonia production facilities.

Overall, this work provides a solid foundation for understanding and improving one of the most important industrial chemical processes, contributing to enhanced efficiency and sustainability in global fertilizer and feedstock production.